\documentclass[graybox]{svmult}

\usepackage{mathptmx}       
\usepackage{helvet}         
\usepackage{courier}        
\usepackage{type1cm}        
\usepackage{makeidx}         
\usepackage{graphicx}        
\usepackage{multicol}        
\usepackage[bottom]{footmisc}

\usepackage{amsfonts}
\usepackage{mathrsfs}
\usepackage{amssymb}
\usepackage{latexsym}
\usepackage{graphicx}
\usepackage{enumerate}
\usepackage{amsmath}
\usepackage{mathtools}
\usepackage{color}
\usepackage{float}
\usepackage{array}
\usepackage{stmaryrd}

\usepackage{eucal}
\usepackage{wrapfig}

\usepackage[round]{natbib}

\usepackage{amssymb,tikz}

\newcommand{\defin}[1]{\textbf{#1}}

\newcommand{\lthen}{\rightarrow}
\newcommand{\liff}{\leftrightarrow}
\newcommand{\falsum}{\bot}
\newcommand{\verum}{\top}
\newcommand{\proves}{\vdash}
\newcommand{\defeq}{\coloneqq}

\newcommand{\val}[1]{\llbracket #1 \rrbracket}

\renewcommand{\phi}{\varphi}
\newcommand{\rimp}{\Rightarrow}
\newcommand{\dimp}{\Leftrightarrow}

\newcommand{\commentout}[1]{}

\renewcommand{\oval}[1]{\val{#1}}

\newcommand{\X}{\mathcal{X}}
\newcommand{\Y}{\mathcal{Y}}
\newcommand{\B}{\mathcal{B}}
\renewcommand{\S}{\mathcal{S}}
\newcommand{\M}{\mathcal{M}}
\newcommand{\T}{\mathcal{T}}
\newcommand{\C}{\mathcal{C}}
\renewcommand{\int}{\mathit{int}}
\newcommand{\spre}{\mathsf{pre}}
\newcommand{\sint}{\mathsf{int}}
\newcommand{\bisim}{\rightleftharpoons}

\makeindex             

\begin{document}

\title*{Topological Subset Space Models for Public Announcements}
\author{Adam Bjorndahl}
\institute{Adam Bjorndahl \at Carnegie Mellon University, 5000 Forbes Avenue Pittsburgh, PA 15213.\\ \email{abjorn@andrew.cmu.edu}
}
%
%
\maketitle


\abstract{We reformulate a key definition given by \citet{WA13} to provide semantics for public announcements in subset spaces. More precisely, we interpret the \emph{precondition} for a public announcement of $\phi$ to be the ``local truth'' of $\phi$, semantically rendered via an interior operator. This is closely related to the notion of $\phi$ being ``knowable''. We argue that these revised semantics improve on the original and offer several motivating examples to this effect. A key insight that emerges is the crucial role of topological structure in this setting. Finally, we provide a simple axiomatization of the resulting logic and prove completeness.}


\section{Introduction} \label{sec:int}

In the standard semantics for epistemic logic, knowledge is represented in terms of possibility: associated with each world $w$ is a set of worlds $R(w)$ representing those states of affairs that are compatible with the agent's knowledge; the agent is said to \emph{know} $\phi$ at $w$ just in case $\phi$ is true at all worlds in $R(w)$ \citep{Hintikka62}. In this context, a world represents a particular arrangement of facts, while a set of worlds represents a particular state of knowledge.

Consider now a \textit{set of sets} of worlds, $\S$: such an object might be construed as representing not how things \textit{are} or what is \textit{known}, but what is \textit{knowable}. Roughly speaking, by restricting attention to models in which each $R(w) \in \S$, we constrain the possible states of knowledge to exactly those in $\S$.

Subset space semantics \citep{DMP96} put this intuition at center stage. In this formalism, the usual relation $R$ is replaced with a collection $\S$ as above, and formulas are evaluated with respect to world-set pairs $(w,U)$, where $w \in U \in \S$, rather than just worlds. In this context, $U$ is called the \emph{epistemic range}. Thus, the possible states of knowledge become an explicit parameter of the model. This provides a convenient setting for studying the \emph{dynamics} of knowledge: learning something new can be captured by shrinking the epistemic range, for instance by transitioning from $(w,U)$ to $(w,V)$, where $V \subseteq U$. Such dynamics are a core concern of subset space logic, which includes an \emph{epistemic effort} modality quantifying over all ways of shrinking the epistemic range in order to express this abstract notion of learning.

One concrete and popular manifestation of epistemic effort is that which results from a \textit{public announcement} \citep{Plaza07}. Intuitively, in the case of a single agent, a public announcement of $\phi$ simply causes the agent to learn that $\phi$ is (or was) true. Subset space models, being well-suited to implementing epistemic updates as discussed above, are a natural and appealing framework in which to interpret public announcements. Somewhat surprisingly, it is only quite recently that this project has been been taken up. \citet{BvDK13} interpret an announcement of $\phi$ in subset space models by essentially the same mechanism as in more standard settings: namely, by deleting those objects that do not satisfy $\phi$. \citet{WA13}, by contrast, interpret announcements using the learning mechanism that is built into the definition of subset space models: that is, by shrinking the epistemic range. It is this latter approach that we focus on.

This article presents a \emph{topological} reformulation of the semantics for public announcements given by \citet{WA13}. The presence of topological structure is convenient for many applications of interest, but the motivation for this project runs deeper: I argue that topology is an \textit{essential} ingredient for the appropriate interpretation of public announcements. This argument is based on two related criticisms of the model presented by \citeauthor{WA13}. First, the \emph{preconditions} they impose for announcements are too strong: certain formulas that really ought to be announceable in their system are not (see Example \ref{exa:tw2}). Second, the epistemic updates produced by successful announcements are not strong enough: loosely speaking, one ought to be able to infer from an announcement of $p$ not only that $p$ is true, but that $p$ is \textit{knowable} (see Example \ref{exa:jat}). In both cases, the resolution of these criticisms motivates and relies upon the foundational topological notion of ``local truth''.

The rest of the paper is organized as follows. In Section \ref{sec:pre}, I present the basics of subset spaces and public announcements, and review the semantics defined by \citet{WA13}. In Section \ref{sec:tsm}, I motivate a revision to these semantics by exhibiting some key interpretational difficulties they face; I then define a topological reformulation and show how it resolves these issues. Section \ref{sec:tch} presents technical results associated with the new topological semantics, including a sound and complete axiomatization. Section \ref{sec:dis} concludes with a discussion of related and future work.


\section{Preliminaries} \label{sec:pre}

\subsection{Subset Space Semantics} \label{sec:sss}

A \defin{subset space} is a pair $(X,\S)$ where $X \neq \emptyset$ is a set of \emph{worlds} (or \emph{states}, or \emph{points}, etc.) and $\S \subseteq 2^{X}$ is a collection of subsets of $X$. Intuitively, elements of $X$ represent \emph{ways the world might be}, while sets in $\S$ represent possible \emph{states of knowledge}.

To make these intuitions precise, consider the (single-agent) epistemic language $\mathcal{EL}$ recursively generated by the grammar
$$
\phi ::= p \, | \, \lnot \phi \, | \, \phi \land \psi \, | \, K \phi,
$$
where $p \in \textsc{prop}$, the (countable) set of \emph{primitive propositions}. Read $K \phi$ as ``the agent knows $\phi$''. A \defin{subset model} $\X = (X,\S,v)$ is a subset space $(X,\S)$ together with a function $v: \textsc{prop} \to 2^{X}$ specifying, for each primitive proposition $p \in \textsc{prop}$, its \emph{extension} $v(p)$. Truth is evaluated with respect to \emph{epistemic scenarios}, which are pairs of the form $(x,U)$, where $x \in U \in \S$. Let $ES(\X)$ denote the collection of all such pairs in $\X$. Given an epistemic scenario $(x,U) \in ES(\X)$, the set $U$ is called its \emph{epistemic range}; it functions like an information set in the sense that knowledge statements at $(x,U)$ are evaluated by universal quantification over $U$. More precisely, we interpret $\mathcal{EL}$ in $\X$ as follows:\footnote{The original definition of subset models \citep{DMP96} was largely motivated by their use in interpreting a richer language containing a second modality representing ``epistemic effort''; roughly speaking, this modality works by shrinking the epistemic range. In the present context, following \citet{WA13} and in the spirit of \citet{BBvDHHdL08}, the mantle of ``epistemic effort'' is taken up by public announcements, so we omit the abstract effort modality. We return to discuss this further in Section \ref{sec:dis}.}
$$
\begin{array}{lcl}
(\X,x,U) \models p & \textrm{ iff } & x \in v(p)\\
(\X,x,U) \models \lnot \phi & \textrm{ iff } & (\X,x,U) \not\models \phi\\
(\X,x,U) \models \phi \land \psi & \textrm{ iff } & (\X,x,U) \models \phi \textrm{ and } (\X,x,U) \models \psi\\
(\X,x,U) \models K \phi & \textrm{ iff } & (\forall y \in U)((\X,y,U) \models \phi).
\end{array}
$$
We sometimes drop mention of the subset model $\X$ when it is clear from context. To get a better sense of how subset space semantics work, an example is helpful.

\begin{example}[The Target and the Wall] \label{exa:tw1}
Consider Figure \ref{fgr:taw},
\begin{figure}
\includegraphics[width=\textwidth]{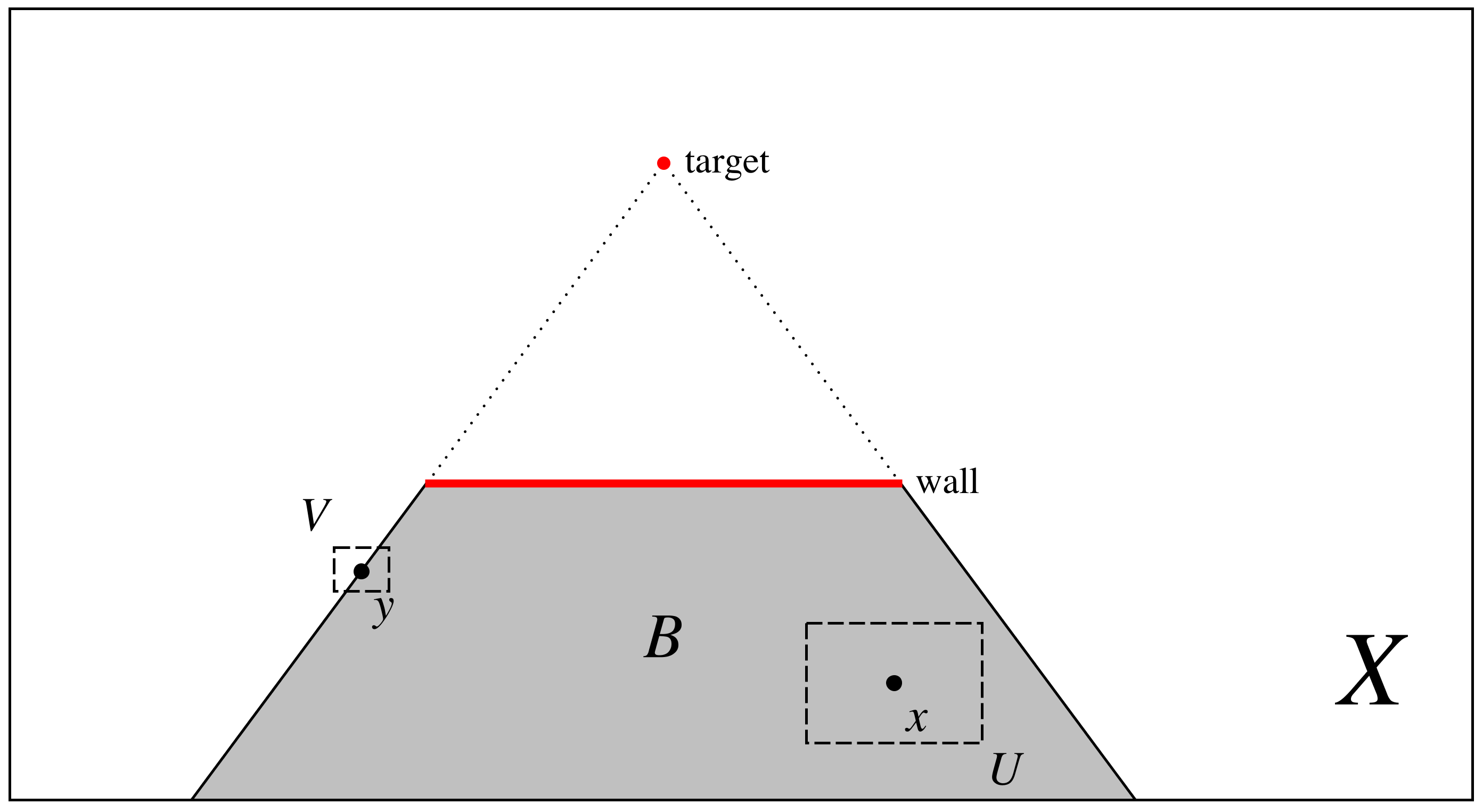}
\caption{A rectangular room with a target and a wall} \label{fgr:taw}
\end{figure}
depicting a rectangular room $X$ into which you have launched a probe. You don't know exactly where it landed, but the probe can measure its distance from the sides of the room and send this data back to you. Of course, any such measurements come with some error. For example, though the probe may have landed at the point $x$, its measurements might only indicate that it is between $0.5$ and $1.5$ meters from the south wall, and between $4$ and $5.5$ meters from the east wall. This can be represented with a rectangle $U$, as shown. 

Let $\S$ consist of those regions of $X$ that can be picked out as above; in other words, if we think of $X$ as a region in $\mathbb{R}^{2}$, $\S$ can be defined as the set of all rectangles $(a,b) \times (c,d) \subseteq X$, where $a < b$ and $c < d$. This definition allows us to put subset space semantics to work in formalizing our intuitions about the knowledge and uncertainty of the person who launched the probe. Suppose, for example, that the probe landed at the point $x$ and returned the measurements above: this corresponds to the epistemic scenario $(x,U)$. On the other hand, if the probe landed at $x$ but failed to return any measurements, this would correspond to the epistemic scenario $(x,X)$. The sets $U$ and $X$ represent the uncertainty that results from different measurements.

There is a target in the room as well as a wall. Assume that you know the location of these objects in advance. From certain vantage points within the room, the wall blocks the target; the shaded region $B$ denotes the set of points where this is so. We might then think of $B$ as the extension of a primitive proposition $b \in \textsc{prop}$ that says ``the wall is blocking the target''.

The relationship between the measurements returned by the probe and your state of knowledge regarding whether the wall is blocking the target is borne out by the semantics defined above. In particular, in the scenario where you receive measurements implying that the probe is in the region $U$, you ought to \textit{know} on the basis of these measurements that the wall is blocking the target, and indeed we have $(x,U) \models Kb$. By contrast, in the scenario where you receive no measurements at all, intuitively, you do not know whether the wall is blocking the target, and this corresponds to the fact that $(x,X) \models \lnot Kb \land \lnot K \lnot b$. Similarly, if the probe lands at the point $y \in B$ (where, intuitively, the wall is ``just barely'' blocking the target) and returns measurements indicating it lies in the region $V$, in the corresponding epistemic scenario $(y,V)$ we also have $(y,V) \models \lnot Kb \land \lnot K \lnot b$. In fact, in this case, since $y$ lies on the borderline between $B$ and its complement, we can see that no measurement, no matter how precise, will yield knowledge of $b$ or its negation. \qed
\end{example}

\subsection{Interpreting Public Announcements in Subset Models} \label{sec:ipa}

We next review the basics of public announcements and the semantics offered by \citet{WA13} for interpreting them in subset models. The (single-agent) public announcement language, denoted $\mathcal{PAL}$, is recursively generated by the grammar
$$
\phi ::= p \, | \, \lnot \phi \, | \, \phi \land \psi \, | \, K \phi \, | \, [\phi] \psi,
$$
where $p \in \textsc{prop}$. The formula $[\phi] \psi$ is read, ``after an announcement of $\phi$, $\psi$ (is true)''. Traditionally, the interpretation of this formula is of the general form
$$
(\M, \omega) \models [\phi]\psi \quad \textrm{iff} \quad (\M, \omega) \models \phi \, \rimp \, (\M|_{\phi}, \omega) \models \psi,
$$
where, loosely speaking, $\M|_{\phi}$ denotes the model obtained from $\M$ by deleting those truth-bearing objects (e.g., worlds) in $\M$ that do not satisfy $\phi$ \citep[see, e.g.,][Chapter 4]{vDvdHK08}. In other words, provided $\phi$ is true, $[\phi]\psi$ holds just in case $\psi$ is true when all $\lnot \phi$ possibilities are removed.

As we have observed, subset spaces offer a model-internal mechanism for representing states of knowledge that obtain ``after some effort''. \citet{WA13} leverage this fact to define an interpretation of public announcements in subset models that implements the update by shrinking the epistemic range rather than by altering the model itself. As a first attempt at defining such a semantics, we might consider the following:

\begin{equation} \label{eqn:fst}
(\X,x,U) \models [\phi]\psi \quad \textrm{iff} \quad (\X,x,U) \models \phi \, \rimp \, (\X,x,\oval{\phi}^{U}) \models \psi,
\end{equation}
where
$$\oval{\phi}^{U} \defeq \{y \in U \: : \: (\X,y,U) \models \phi\},$$
called the \defin{extension of $\phi$ under $U$}. The idea is that shrinking the epistemic range from $U$ to $\oval{\phi}^{U}$ captures the effect of hearing a public announcement of $\phi$. An immediate problem with this definition is that $(x,\oval{\phi}^{U})$ may not be an epistemic scenario: it is if and only if $\oval{\phi}^{U} \in \S$.

Call the antecedent of the implication in (\ref{eqn:fst}) the \emph{precondition} for the announcement, and the consequent the \emph{postcondition}. The definition proposed by \citeauthor{WA13} avoids the issue raised above by strengthening the precondition in such a way as to ensure that the postcondition is defined:
\begin{equation} \label{eqn:wad}
(\X,x,U) \models [\phi]\psi \quad \textrm{iff} \quad (\X,x,U) \models \spre(\phi) \, \rimp \, (\X,x,\oval{\phi}^{U}) \models \psi,
\end{equation}
where
$$(\X,x,U) \models \spre(\phi) \quad \textrm{iff} \quad x \in \oval{\phi}^{U} \in \S.$$

These semantics offer a way of interpreting public announcements that obviates the need to consider alternative models, thus elegantly realizing the central insight of \citeauthor{WA13}. Clearly, $\spre(\phi)$ strengthens the classical precondition, which simply insists that $\phi$ be true. However, as we now show, this precondition is in fact \emph{too} strong, and moreover, the postcondition is too weak.

\section{Topological Subset Models} \label{sec:tsm}

\subsection{Motivation}

For $\spre(\phi)$ to hold in an epistemic scenario $(x,U)$, two conditions must be satisfied. First, $x$ must be in $\oval{\phi}^{U}$, which is simply the subset model analogue of the classical precondition that $\phi$ be true. Second, we must have $\oval{\phi}^{U} \in \S$. It is tempting to read this latter condition as something like, ``$\phi$ is knowable (given $U$)''. After all, $\S$ collects precisely those subsets of $X$ that can function as states of knowledge. But this reading is misleading: as the examples below make clear, it is possible to know $\phi$ even if its extension is not a member of $\S$. Speaking abstractly, to insist that $\oval{\phi}^{U} \in \S$ is to impose a ``global'' precondition on announcements where we should instead be appealing to a ``local'' condition. To make these ideas concrete, we return to the setting of Example \ref{exa:tw1}.

\begin{example}[The Target and the Wall, continued] \label{exa:tw2}
Recall that you have launched a probe into a room containing a target and a wall, as depicted in Figure \ref{fgr:taw}. Due to our definition of $\S$, which effectively identifies states of knowledge with certain types of measurements, we have $\oval{b}^{X} = B \notin \S$ (since $B$ is not a rectangle). This implies that for all $z \in X$, $(z,X) \not\models \spre(b)$, and so by the definition given in (\ref{eqn:wad}), ``the wall is blocking the target'' is not announceable in any epistemic scenario of the form $(z,X)$.\footnote{More precisely, it means that every formula of the form $[b]\phi$ is trivially true at $(z,X)$, even when $\phi$ is a contradiction.}

This seems wrong: there are some epistemic scenarios of the form $(z,X)$ in which $b$ really ought to be announceable. For instance, suppose that the probe in fact landed at the point $x$ but you have received no measurements, corresponding to the epistemic scenario $(x,X)$. In this case, not only is the wall blocking the target, but this fact is ``knowable'' in the sense that there is a measurement---for example, the rectangle $U$ depicted in Figure \ref{fgr:taw}---that entails it. We might even imagine that some third party has intercepted the probe's transmission of the measurement $U$. An adequate theory of public announcements should predict that this third party can meaningfully announce to you, ``The wall is blocking the target''.

This highlights the ``global versus local'' distinction alluded to above: informally, although $\oval{b}^{X}$ is not itself in $\S$, there are elements of $\S$ that entail $b$, and in the right epistemic scenarios this seems sufficient to license the public announcement. Roughly speaking, we might say that $b$ is \emph{locally true} at $(x,X)$ because $x \in U \subseteq \oval{b}^{X}$, and in general redefine the precondition for an announcement of $\phi$ so that it demands only local truth. We make these notions precise in Section \ref{sec:frm}. \qed
\end{example}

Before turning to the formalism, we consider one more example of a rather different character.

\begin{example}[The Jewel and the Tomb] \label{exa:jat}
You have learned from ancient historical records of the existence of a secret tomb within which was supposedly ensconced a priceless jewel. In point of fact, you have no idea whether a priceless jewel was actually placed within this tomb before it was sealed---perhaps that part of the historical record was simply an embellishment. You are also unsure as to whether this tomb is still lost or has been rediscovered in modern times (and its contents catalogued).

The relevant possibilities here can be captured with a four-state model: let $X = \{s_{JD}, s_{J\bar{D}}, s_{\bar{J}D}, s_{\bar{J}\bar{D}}\}$, where each state in $X$ encodes whether the tomb actually contains a jewel ($J$) or not ($\bar{J}$), and whether it has been rediscovered in modern times ($D$) or not ($\bar{D}$).

We also want our model to encode the fact that \textit{the only way to learn about the jewel is to discover the tomb} (all other records of the jewel's existence, or lack thereof, having been irrevocably lost to time). Subset spaces are ideally suited to encoding such constraints on the possible states of knowledge; this is accomplished by controlling the elements of $\S$. In this example, you could conceivably know whether or not the tomb has been discovered in modern times without knowing whether or not there is a jewel inside, corresponding to the two knowledge states $\{s_{JD}, s_{\bar{J}D}\}$ and $\{s_{J\bar{D}}, s_{\bar{J}\bar{D}}\}$. Furthermore, \textit{provided} you know that the tomb has been discovered, you might also know whether or not a priceless jewel was found inside, corresponding to the two knowledge states $\{s_{JD}\}$ and $\{s_{\bar{J}D}\}$. We therefore define
$$\S = \{\{s_{JD}, s_{\bar{J}D}\}, \{s_{J\bar{D}}, s_{\bar{J}\bar{D}}\}, \{s_{JD}\}, \{s_{\bar{J}D}\}\}.$$
Crucially, we do not have, for example, $\{s_{J\bar{D}}\} \in \S$, since this would correspond to a state of knowledge where you know both that the tomb has not been rediscovered in modern times and that there is a jewel inside. This is precisely what we want to rule out.

Let $j$ and $d$ be primitive propositions standing for ``the jewel is in the tomb'' and ``the tomb has been discovered'', respectively, and let $v: \{j,d\} \to 2^{X}$ be defined in the obvious way. Then it is easy to see that for all $x \in X$, $(x,X) \not\models \spre(\lnot j \land \lnot d)$; this follows from the fact that $\oval{\lnot j \land \lnot d}^{X} = \{s_{\bar{J}\bar{D}}\} \notin \S$. This accords with the intuition that since you cannot know that the jewel is not in the tomb without also knowing that the tomb has been discovered, such a state of affairs should not be announceable.

However, it is also easy to see that $(s_{JD}, X) \not\models \spre(j)$; this follows from the fact that $\oval{j}^{X} = \{s_{JD}, s_{J\bar{D}}\} \notin \S$. Thus, ``the jewel is in the tomb'' is not announceable even if the jewel really is in the tomb \textit{and} the tomb has been discovered. This seems wrong---for instance, the person who discovered the tomb could have seen the jewel inside and then announced this fact. As in the previous example, this intuition is borne out in the notion of local truth: there is a state of knowledge $\{s_{JD}\} \in \S$ that entails $j$, and this ought to be a sufficient condition for the announceability of $j$ in the epistemic scenario $(s_{JD},X)$.

Note also that $\{s_{JD}\}$ is the \textit{only} element of $\S$ that entails $j$; this corresponds to the fact that the only way to learn about the jewel is to discover the tomb. As a consequence, any successful announcement of $j$ should carry with it the \textit{implication} that the tomb has already been discovered. In other words, we ought to have $(s_{JD}, X) \models \spre(j) \land [j]Kd$.

The semantics defined by \citeauthor{WA13} do not give credence to these intuitions. First, as we saw, $\spre(j)$ is not satisfied at $(s_{JD},X)$. Furthermore, since announcements in this framework have the effect of restricting the epistemic range to the extension of the announced formula, no inferences beyond the truth of that formula (and the logical consequences thereof) are supported. In particular, since $\oval{j}^{X} \not\subseteq \oval{d}^{X}$, $j$ does not entail $d$, so even if $j$ were announceable it would not result in $d$ becoming known. This suggests that in addition to weakening the precondition for a public announcement, we also need to strengthen the postcondition. \qed
\end{example}

\subsection{Formal Semantics} \label{sec:frm}

The notion of ``local truth'' is naturally and succinctly captured in a topological framework. A \defin{topological space} is a pair $\X = (X,\T)$ where $X$ is a nonempty set and $\T \subseteq 2^{X}$ is a collection of subsets of $X$ that covers $X$ and is closed under finite intersections and arbitrary unions. The collection $\T$ is called a \emph{topology on $X$} and elements of $\T$ are called \emph{open} sets.

Topology might be described as the abstract mathematics of space. Roughly speaking, each open set can be viewed as encoding a notion of ``nearness''; this notion is operationalized in the following definition. Given a set $A \subseteq X$, we say that $x$ lies in the \defin{interior} of $A$ if there is some $U \in \T$ such that $x \in U \subseteq A$. The open set $U$ acts a ``witness'' to $x$'s membership in $A$: not only is $x$ in $A$, but also all ``nearby'' points (i.e., all $y \in U$) are in $A$.

The set of all points in the interior of $A$ is denoted $\int_{\X}(A)$. To ease notational clutter, we often drop the subscript and sometimes omit the parentheses. It is not hard to see that $\int(A) \in \T$: for each $x \in \int(A)$, there is by definition an open set $U_{x}$ such that $x \in U_{x} \subseteq A$, and it is easy to check that $\bigcup_{x \in \int(A)} U_{x} = \int(A)$. In fact, $\int(A)$ is the largest open set contained in $A$. For a general introduction to topology we refer the reader to \citet{Munkres}.

A \defin{topological subset model} is a subset model $\X = (X,\T,v)$ in which $\T$ is a topology on $X$. Since every topological space is a subset space, the epistemic intuitions for subset spaces apply also to topological spaces---we can identify open sets with measurements, or more generally with states of knowledge. But the additional topological structure allows us to go further: in particular, the notion of local truth motivated in Examples \ref{exa:tw2} and \ref{exa:jat} coincides exactly with the definition of topological interior.

The core proposal of this paper is to interpret public announcements on \textit{topological} spaces according to the following reformulated semantics:
\begin{equation} \label{eqn:ref}
(\X,x,U) \models [\phi]\psi \quad \textrm{iff} \quad (\X,x,U) \models \sint(\phi) \, \rimp \, (\X,x,\int\oval{\phi}^{U}) \models \psi,
\end{equation}
where
\begin{equation} \label{eqn:int}
(\X,x,U) \models \sint(\phi) \quad \textrm{iff} \quad x \in \int\oval{\phi}^{U}.
\end{equation}
To distinguish these semantics from those given in (\ref{eqn:wad}), we refer to them as ``$\sint$-semantics'' and ``$\spre$-semantics'', respectively.

Since $\sint$-semantics make use of the interior operator, they are only defined on topological subset models, though of course $\spre$-semantics also make sense in this setting. Comparing the two is instructive; they differ both in the precondition and the postcondition. We first observe that
$$x \in \oval{\phi}^{U} \in \T \, \rimp \, x \in \int\oval{\phi}^{U} \, \rimp \, x \in \oval{\phi}^{U},$$
and neither of the reverse implications holds in general; it follows that $\sint(\phi)$ is a strictly weaker condition than $\spre(\phi)$ and a strictly stronger condition than $\phi$. This, of course, is by design: as we show below, weakening the precondition in this way provides exactly the leeway needed to render the problematic ``unannounceable'' formulas considered in Examples \ref{exa:tw2} and \ref{exa:jat} announceable.

The postcondition, on the other hand, has been strengthened: the updated epistemic range in (\ref{eqn:ref}), namely $\int\oval{\phi}^{U}$, is a subset of the epistemic range in (\ref{eqn:wad}), which is just $\oval{\phi}^{U}$. Note that $\oval{\phi}^{U}$ might not be open in our semantics (even when the precondition is satisfied), so it cannot, in general, serve as an epistemic range. Replacing it with its interior is a convenient fix for this technical issue. But there is a deeper motivation and broader import for this definition: a successful announcement of $\phi$ can carry more information than simply the content of $\phi$ itself. By replacing $\oval{\phi}^{U}$ with $\int\oval{\phi}^{U}$ in the postcondition, we are effectively updating the agent's knowledge with not merely with the \textit{truth} of $\phi$, but with the \textit{announceability} of $\phi$. Said differently: when an agent hears a public announcement of $\phi$, they can deduce not only that $\phi$ is true, but that $\phi$ is entailed by some state of knowledge---they come to know that the true state of the world is somewhere in
$$\bigcup \{V \in \T \: : \: V \subseteq \oval{\phi}^{U}\},$$
which is exactly $\int\oval{\phi}^{U}$.

We explore these features of our semantics in the following examples. A preliminary definition is useful: given a collection of subsets $\C \subseteq 2^{X}$ that covers $X$, the \defin{topology generated by $\C$}, denoted $\T(\C)$, is simply the smallest topology on $X$ containing $\C$. It is not hard to check that $\T(\C)$ is equal to the set of all arbitrary unions of finite intersections of members of $\C$.

\begin{example}[The Target and the Wall, revisited] \label{exa:tw3}
We first transform the subset model given in Example \ref{exa:tw1} into a \textit{topological} model by replacing $\S$ with the topology it generates, $\T(\S)$. Since $\S$ is already closed under finite intersections, this amounts simply to closing under unions. It is not hard to see that $\T(\S)$ is the standard Euclidean topology on the plane relativized to $X$. 

Expanding $\S$ to $\T(\S)$ does not in itself solve the problems raised in Example \ref{exa:tw2}: it is still the case, for instance, that $(x,X) \not\models \spre(b)$ (since $\oval{b}^{X} = B$ is not open). But the presence of topological structure allows us to switch to $\sint$-semantics. Since $x$ lies in the interior of $B$, we have $(x,X) \models \sint(b)$; it follows that ``the wall is blocking the target'' \textit{is} announceable at $(x,X)$ according to $\sint$-semantics, as intuition suggests it ought to be. By contrast, a probe that landed at $y$ would be incapable of transmitting any measurement that entails $b$; as such, we might expect that $b$ is not announceable at $(y,X)$, and indeed, we have $(y,X) \not\models \sint(b)$ (since $y \notin \int(B)$).

Finally, we observe that after a successful announcement of $b$ at $(x,X)$, the updated epistemic range is not the extension of $b$, but rather its interior, $\int(B)$. This corresponds to the idea that an announcement of $b$ in this epistemic scenario conveys more information than just the truth of $b$: it tells you in addition that the probe must have landed at a point where it can actually take some measurement that entails $b$. There are, of course, many such measurements---any open rectangle contained in $B$ is such a measurement. The crucial point is this: to know that \textit{some one} of these measurements must have been taken, but not \textit{which} one in particular, is to know that the true state of the world lies in their union, $\int(B)$. \qed
\end{example}

\begin{example}[The Jewel and the Tomb, revisited] \label{exa:jt2}
Once again, we extend the collection $\S$ given in Example \ref{exa:jat} to the topology it generates:
$$\T(\{\{s_{JD}, s_{\bar{J}D}\}, \{s_{J\bar{D}}, s_{\bar{J}\bar{D}}\}, \{s_{JD}\}, \{s_{\bar{J}D}\}\}).$$
As with $\spre$-semantics, $\sint$-semantics determines that $\lnot j \land \lnot d$ is not announceable in any epistemic scenario of the form $(x,X)$; this follows from the fact that $\int(\{s_{\bar{J}\bar{D}}\}) = \emptyset$. By contrast, ``the jewel is in the tomb'' \textit{is} announceable in the epistemic scenario $(s_{JD}, X)$---as it ought to be---on account of the fact that
$$\int\oval{j}^{X} = \int(\{s_{JD}, s_{J\bar{D}}\}) = \{s_{JD}\}.$$
Observe also that the epistemic range in the postcondition for announcing $j$ is $\int\oval{j}^{X} = \{s_{JD}\}$. Since $(s_{JD},\{s_{JD}\}) \models Kd$, we therefore have $(s_{JD}, X) \models [j]Kd$, which is exactly the effect we sought in Example \ref{exa:jat}: the only way to learn about the jewel is to discover the tomb, so any successful announcement of $j$ should carry with it the implication that the tomb has indeed been discovered. \qed
\end{example}

\section{Technical Results} \label{sec:tch}

Throughout this section, except where otherwise noted, we work with $\sint$-semantics as given in (\ref{eqn:ref}). We say that $\phi$ is \defin{valid} and write $\models \phi$ if, for all topological subset models $\X$ and all epistemic scenarios $(x,U) \in ES(\X)$, we have $(\X, x, U) \models \phi$.

We begin by establishing some basic properties of the $\sint$ modality. Observe first that $\sint(\phi)$ is definable in $\mathcal{PAL}$: it is semantically equivalent to the formula $\lnot [\phi] \falsum$ (where $\falsum$ denotes some propositional contradiction). As such, we can freely add the $\sint$ modality to $\mathcal{PAL}$ without changing the expressivity of the language, and in the following we take this for granted.

The properties of a modalized interior operator have been thoroughly investigated \citep[see, e.g.,][]{AvBB03,vBB07}, so much of the following proposition should come as little surprise.
\begin{proposition} \label{pro:prp}
For all $\phi, \psi \in \mathcal{PAL}$, the following hold:
\begin{enumerate}[(a)]
\item
$\models \sint(\phi) \lthen \phi$
\item
$\models \sint(\phi) \lthen \sint(\sint(\phi))$
\item
$\models \sint(\phi \lthen \psi) \lthen (\sint(\phi) \lthen \sint(\psi))$
\item
$\models \phi$ implies $\models \sint(\phi)$
\item
$\models K \phi \lthen \sint(\phi)$
\item
$\not\models \sint(\phi) \lthen K(\phi \lthen \sint(\phi))$
\item
$\not\models \lnot(\phi \lthen \sint(\phi)) \lthen K \lnot \sint(\phi)$.
\end{enumerate}
\end{proposition}

\begin{proof}
Parts (a) through (d) constitute a standard \textsf{S4} axiomatization of the interior operator, and the proof that they hold in this setting is analogous to the usual proof. Part (e) follows from the fact that
\begin{eqnarray*}
(\X,x,U) \models K \phi & \, \rimp \, & \oval{\phi}^{U} = U\\
& \, \rimp \, & int\oval{\phi}^{U} = U.
\end{eqnarray*}
\begin{minipage}{.65\textwidth}
\vspace{-6mm}
Parts (f) and (g) are included to exhibit some of the differences between $\spre(\phi)$ (for which these two schemes are valid) and $\sint(\phi)$. Let $X$ be a subset of the plane equipped with the standard Euclidean subspace topology. Assume that $U$ is an open subset of $X$, as shown in Figure \ref{fgr:cex}, and let $v(p) = U \cup \{y\}$. Then it is easy to check that $(x, X) \not \models \sint(p) \lthen K(p \lthen \sint(p))$ and $(y,X) \not\models \lnot(p \lthen \sint(p)) \lthen K \lnot \sint(p)$. \qed
\end{minipage}
\hfill
\begin{minipage}{.3\textwidth}
\vspace{-6mm}
\begin{figure}[H]
\includegraphics[width=\textwidth]{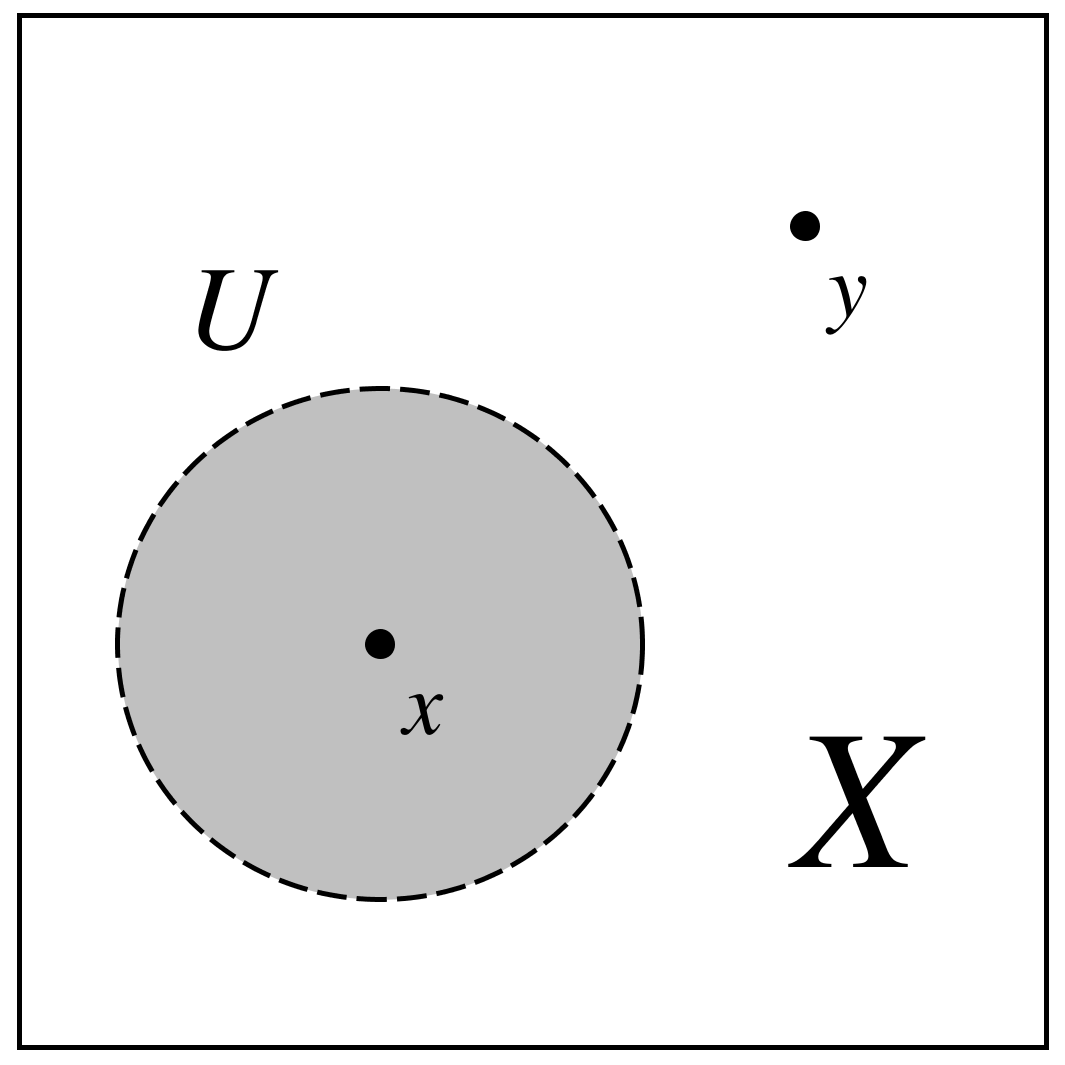}
\caption{Counterexamples} \label{fgr:cex}
\end{figure}
\end{minipage}
\end{proof}


Although $\sint$ is definable in $\mathcal{PAL}$, it also makes sense to consider in a language without public announcements. This plays an important role in our axiomatization. Let $\mathcal{EL}_{\sint}$ be recursively generated by the grammar
$$\phi ::= p \, | \, \lnot \phi \, | \, \phi \land \psi \, | \, K \phi \, | \, \sint(\phi),$$
where $p \in \textsc{prop}$. This language is interpreted in topological subset models in the obvious way; in particular, the semantics of $\sint$ are given as before by (\ref{eqn:int}).

$\mathcal{EL}_{\sint}$ is an extension of $\mathcal{EL}$ and is strictly more expressive, since $\mathcal{EL}$ cannot define $\sint(\phi)$. To show this, we first recall the following definition \citep[Definition 8]{WA13}: given two subset models $\X$ and $\X'$, a relation $\bisim$ between $ES(\X)$ and $ES(\X')$ is called a \defin{partial bisimulation} (between $\X$ and $\X'$) if whenever $(x,U) \bisim (x',U')$, the following conditions are satisfied:
\begin{description}[\textsc{Forth}]
\item[\textsc{Base}]
$(\forall p \in \textsc{prop})(x \in v(p) \dimp x' \in v'(p))$
\item[\textsc{Forth}]
$(\forall y \in U)(\exists y' \in U')((y,U) \bisim (y',U'))$
\item[\textsc{Back}]
$(\forall y' \in U')(\exists y \in U)((y,U) \bisim (y',U'))$.
\end{description}
This is the natural analogue of the usual notion of bisimulation defined on relational structures \citep[see, e.g.,][]{BRV01}. An easy structural induction over $\mathcal{EL}$ yields the following invariance result.

\begin{proposition} \label{pro:bis}
Let $\bisim$ be a partial bisimulation between subset models $\X$ and $\X'$ with $(x,U) \bisim (x',U')$. Then for all $\phi \in \mathcal{EL}$,
$$(\X, x, U) \models \phi \, \dimp \, (\X', x', U') \models \phi.$$
\end{proposition}

\begin{proposition}
$\mathcal{EL}_{\sint}$ is strictly more expressive than $\mathcal{EL}$.
\end{proposition}

\begin{proof}
By Proposition \ref{pro:bis}, it suffices to show that $\sint(p)$ can distinguish two epistemic scenarios that are linked by a partial bisimulation. Consider the topological subset models
$$\X = (\{x,y\}, 2^{\{x,y\}}, v)$$
and
$$\Y = (\{x,y\}, \{\emptyset, \{y\}, \{x,y\}\}, v),$$
where $v(p) = \{x\}$. Thus, in $\X$ we have the discrete topology, while in $\Y$ the singleton $\{y\}$ is open but $\{x\}$ is not. It is easy to check that the relation given by
\begin{eqnarray*}
(x, \{x,y\}) & \bisim & (x, \{x,y\})\\
(y, \{x,y\}) & \bisim & (y, \{x,y\})
\end{eqnarray*}
is a partial bisimulation. However, we have $(\X, x, \{x,y\}) \models \sint(p)$, since
$$x \in \{x\} = \int_{\X}(\{x\}) = \int_{\X}\oval{p}^{\{x,y\}},$$
whereas $(\Y, x, \{x,y\}) \not\models \sint(p)$, since
$$x \notin \emptyset = \int_{\Y}(\{x\}) = \int_{\Y}\oval{p}^{\{x,y\}}.$$
Hence, $\sint(p)$ cannot be equivalent to any formula of $\mathcal{EL}$.
\qed
\end{proof}

Of course, this result also shows that $\mathcal{PAL}$ is strictly more expressive than $\mathcal{EL}$ (in $\sint$-semantics). On the other hand, $\mathcal{EL}_{\sint}$ and $\mathcal{PAL}$ are \textit{equally} expressive: in essence, this is because the following reduction schemes allow us to rewrite any formula of $\mathcal{PAL}$ as a logically equivalent formula of $\mathcal{EL}_{\sint}$ \citep[cf.][Theorem 11]{WA13}.

\begin{proposition} \label{pro:red}
The following $\mathcal{PAL}$ formulas are valid:
$$
\begin{array}{rcl}
{[\phi] p} & \liff & (\sint(\phi) \lthen p)\\
{[\phi] \lnot \psi} & \liff & (\sint(\phi) \lthen \lnot {[\phi] \psi})\\
{[\phi] (\psi \land \chi)} & \liff & ({[\phi] \psi} \land {[\phi] \chi})
\end{array}
\quad
\begin{array}{rcl}
{[\phi] K \psi} & \liff & (\sint(\phi) \lthen K {[\phi] \psi})\\
{[\phi] \sint(\psi)} & \liff & (\sint(\phi) \lthen \sint({[\phi] \psi}))\\
{[\phi][\psi] \chi} & \liff & [\sint(\phi) \land {[\phi] \sint(\psi)}] \chi.
\end{array}
$$
\end{proposition}

\begin{proof}
The first three equivalences are straightforward to prove. To show that ${[\phi] K \psi} \liff (\sint(\phi) \lthen K {[\phi] \psi})$ is valid, first note that if $(x,U) \not\models \sint(\phi)$ then this equivalence holds trivially at $(x,U)$. Otherwise, assuming that $(x,U) \models \sint(\phi)$, we have:
\begin{eqnarray*}
(x,U) \models {[\phi] K \psi} & \dimp & (x,\int\oval{\phi}^{U}) \models K \psi\\
& \dimp & (\forall y \in \int\oval{\phi}^{U})((y, \int\oval{\phi}^{U}) \models \psi),
\end{eqnarray*}
whereas
\begin{eqnarray*}
(x,U) \models \sint(\phi) \lthen K {[\phi] \psi} & \dimp & (x,U) \models K {[\phi] \psi}\\
& \dimp & (\forall z \in U)((z,U) \models {[\phi] \psi})\\
& \dimp & (\forall z \in U)(z \in \int\oval{\phi}^{U} \rimp (z,\int\oval{\phi}^{U}) \models \psi)\\
& \dimp & (\forall z \in \int\oval{\phi}^{U})((z,\int\oval{\phi}^{U}) \models \psi).
\end{eqnarray*}

Next we show that $[\phi] \sint(\psi) \liff (\sint(\phi) \lthen \sint([\phi] \psi))$ is valid. As above, this equivalence holds trivially at $(x,U)$ when $(x,U) \not\models \sint(\phi)$, so assume that $(x,U) \models \sint(\phi)$. We then have:
\begin{eqnarray}
(x,U) \models {[\phi] \sint(\psi)} & \dimp & (x,\int\oval{\phi}^{U}) \models \sint(\psi) \nonumber \\
& \dimp & x \in \int\oval{\psi}^{\int\oval{\phi}^{U}} \nonumber \\
& \dimp & x \in \int\{y \in \int\oval{\phi}^{U} \: : \: (y,\int\oval{\phi}^{U}) \models \psi\}, \label{eqn:in1}
\end{eqnarray} 
and
\begin{eqnarray}
(x,U) \models \sint(\phi) \lthen \sint({[\phi] \psi}) & \dimp & (x,U) \models \sint({[\phi] \psi}) \nonumber\\
& \dimp & x \in \int\oval{[\phi]\psi}^{U} \nonumber\\
& \dimp & (\exists V \in \T)(x \in V \subseteq \oval{[\phi] \psi}^{U}). \label{eqn:in2}
\end{eqnarray}
Now observe that
\begin{eqnarray*}
\oval{[\phi] \psi}^{U} & = & \{y \in U \: : \: (y,U) \models [\phi] \psi\}\\
& = & \{y \in U \: : \: y \in \int\oval{\phi}^{U} \rimp (y,\int\oval{\phi}^{U}) \models \psi\},
\end{eqnarray*}
so clearly any witness $V \in \T$ to (\ref{eqn:in1}) also satisfies (\ref{eqn:in2}). Conversely, given a $V$ satisfying (\ref{eqn:in2}), let $V' = V \cap \int\oval{\phi}^{U}$. By assumption, $x \in \int\oval{\phi}^{U}$, so we have $x \in V'$, and it is easy to see that $V'$ is a witness to (\ref{eqn:in1}).

Finally, to see that ${[\phi][\psi] \chi} \liff [\sint(\phi) \land {[\phi] \sint(\psi)}] \chi$ is valid, first observe that
\begin{eqnarray*}
(x,U) \models {[\phi][\psi] \chi} & \dimp & x \in \int\oval{\phi}^{U} \rimp (x,\int\oval{\phi}^{U}) \models {[\psi] \chi}\\
& \dimp & x \in \int\oval{\phi}^{U} \rimp \big(x \in \int\oval{\psi}^{\int\oval{\phi}^{U}} \rimp (x,\int\oval{\psi}^{\int\oval{\phi}^{U}}) \models \chi\big)\\
& \dimp & x \in \int\oval{\psi}^{\int\oval{\phi}^{U}} \rimp (x,\int\oval{\psi}^{\int\oval{\phi}^{U}}) \models \chi,
\end{eqnarray*}
where the last line follows from the fact that
\begin{equation} \label{eqn:cnt}
\int\oval{\psi}^{\int\oval{\phi}^{U}} \subseteq \oval{\psi}^{\int\oval{\phi}^{U}} \subseteq \int\oval{\phi}^{U}.
\end{equation}
On the other hand, $(x,U) \models {[\sint(\phi) \land {[\phi] \sint(\psi)}] \chi}$ iff
$$x \in \int\oval{\sint(\phi) \land {[\phi] \sint(\psi)}}^{U} \rimp (x,\int\oval{\sint(\phi) \land {[\phi] \sint(\psi)}}^{U}) \models \chi;$$
thus, to complete the proof it suffices to show that
\begin{equation} \label{eqn:sff}
\int\oval{\sint(\phi) \land {[\phi] \sint(\psi)}}^{U} = \int\oval{\psi}^{\int\oval{\phi}^{U}}.
\end{equation}
By definition,
\begin{eqnarray*}
\oval{\sint(\phi) \land {[\phi] \sint(\psi)}}^{U} & = & \{y \in U \: : \: (y,U) \models \sint(\phi) \land {[\phi] \sint(\psi)}\}\\
& = & \{y \in U \: : \: y \in \int\oval{\phi}^{U} \textrm{ and }\\
& & \qquad \qquad (y \in \int\oval{\phi}^{U} \rimp (y,\int\oval{\phi}^{U}) \models \sint(\psi))\}\\
& = & \{y \in \int\oval{\phi}^{U} \: : \: (y,\int\oval{\phi}^{U}) \models \sint(\psi))\}\\
& = & \{y \in \int\oval{\phi}^{U} \: : \: y \in \int\oval{\psi}^{\int\oval{\phi}^{U}}\}\\
& = & \int\oval{\psi}^{\int\oval{\phi}^{U}},
\end{eqnarray*}
where the third line follows from the fact that $\int\oval{\phi}^{U} \subseteq U$, and the last line follows from (\ref{eqn:cnt}). Since $\int^{2} = \int$, this establishes (\ref{eqn:sff}).
\qed
\end{proof}

It remains to show that these reduction schemes actually allow us to rewrite any $\mathcal{PAL}$ formula as an equivalent $\mathcal{EL}_{\sint}$ formula. For this, the following definition is useful \citep[cf.][Definition 22]{WA13}: the \defin{complexity} $c(\phi)$ of any $\mathcal{PAL}$ formula $\phi$ is defined recursively by
$$
\begin{array}{rcl}
{c(p)} & = & 1\\
{c(\lnot \phi)} & = & c(\phi) + 1\\
{c(\phi \land \psi)} & = & c(\phi) + c(\psi) + 1
\end{array}
\quad
\begin{array}{rcl}
{c(K \phi)} & = & c(\phi) + 1\\
{c(\sint(\phi))} & = & c(\phi) + 1\\
{c([\phi] \psi)} & = & (c(\phi) + 6) \cdot c(\psi).
\end{array}
$$
\begin{lemma} \label{lem:com}
Each of the six reduction schemes in Proposition \ref{pro:red} reduces complexity from left to right: the complexity of the formula on the righthand side of the biconditional is less than the complexity of the formula on the lefthand side.
\end{lemma}

\begin{proof}
To begin, observe that
$$c(\phi \lthen \psi) = c(\lnot(\phi \land \lnot \psi)) = c(\phi) + c(\psi) + 3.$$
Now it is easy to check that
$$c(\sint(\phi) \lthen p) = c(\phi) + 5 < c(\phi) + 6 = c([\phi]p).$$
We also have
\begin{eqnarray*}
c(\sint(\phi) \lthen \lnot [\phi]\psi) & = & c(\phi) + 1 + (c(\phi) + 6) \cdot c(\psi) + 1 + 3\\
& = & c(\phi) \cdot c(\psi) + c(\phi) + 6c(\psi) + 5\\
& < & c(\phi) \cdot c(\psi) + c(\phi) + 6c(\psi) + 6\\
& = & (c(\phi) + 6)(c(\psi) + 1)\\
& = & c([\phi] \lnot\psi).
\end{eqnarray*}
The calculations for the reduction schemes corresponding to the $K$ and $\sint$ modalities proceed analogously. Next, we have
\begin{eqnarray*}
c([\phi]\psi \land [\phi]\chi) & = & (c(\phi) + 6) \cdot c(\psi) + (c(\phi) + 6) \cdot c(\chi) + 1\\
& < & (c(\phi) + 6)(c(\psi) + c(\chi) + 1)\\
& = & c([\phi](\psi \land \chi)).
\end{eqnarray*}
And finally:
\begin{eqnarray*}
c([\sint(\phi) \land [\phi]\sint(\psi)]\chi) & = & (c(\phi) + 1 + (c(\phi) + 6)(c(\psi) + 1) + 1 + 6) \cdot c(\chi)\\
& = & (c(\phi) \cdot c(\psi) + 2c(\phi) + 6c(\psi) + 14) \cdot c(\chi)\\
& < & (c(\phi) + 6)(c(\psi) + 6) \cdot c(\chi)\\
& = & c([\phi][\psi]\chi). \quad \qed
\end{eqnarray*}

\end{proof}

\begin{proposition} \label{pro:eqv}
For all $\mathcal{PAL}$ formulas $\phi$, there exists an $\mathcal{EL}_{\sint}$ formula $\tilde{\phi}$ such that $\models \phi \liff \tilde{\phi}$.
\end{proposition}

\begin{proof}
The proof proceeds by induction on $c(\phi)$. If $c(\phi) = 1$, then $\phi \in \textsc{prop}$ and we can take $\tilde{\phi} = \phi$. Now suppose that $c(\phi) > 1$, and assume inductively that the result holds for all formulas with complexity less than $c(\phi)$. There are several cases to consider, depending on the structure of $\phi$.

If $\phi = \lnot \psi$ for some $\psi$, then $c(\psi) < c(\phi)$, so by the inductive hypothesis there is an $\mathcal{EL}_{\sint}$ formula $\tilde{\psi}$ such that $\models \psi \liff \tilde{\psi}$. It follows that $\models \phi \liff \lnot \tilde{\psi}$, which establishes the desired result. The cases corresponding to $\phi = K \psi$ and $\phi = \sint(\psi)$ are handled analogously. The case where $\phi = \psi_{1} \land \psi_{2}$ is also similar: since $c(\psi_{1}) < c(\phi)$ and $c(\psi_{2}) < c(\phi)$, we can find $\mathcal{EL}_{\sint}$ formulas $\tilde{\psi}_{1}$ and $\tilde{\psi}_{2}$ such that $\models \psi_{1} \liff \tilde{\psi}_{1}$ and $\models \psi_{2} \liff \tilde{\psi}_{2}$, hence $\models \phi \liff (\tilde{\psi}_{1} \land \tilde{\psi}_{2})$.

The final case is when $\phi = [\psi]\chi$. By applying one of the reduction schemes in Proposition \ref{pro:red}, we can find a formula $\xi$ such that $\models \phi \liff \xi$; moreover, by Lemma \ref{lem:com}, we know that $c(\xi) < c(\phi)$. The inductive hypothesis now applies to give us an $\mathcal{EL}_{\sint}$ formula $\tilde{\xi}$ such that $\models \xi \liff \tilde{\xi}$. Of course, we then have $\models \phi \liff \tilde{\xi}$, which completes the proof.
\qed
\end{proof}

At last we turn our attention to a sound and complete axiomatization of $\mathcal{PAL}$ in $\sint$-semantics. We first axiomatize $\mathcal{EL}_{\sint}$ and then use the reduction schemes to transform this into an axiomatization of $\mathcal{PAL}$.

Let $\mathsf{CPL}$ denote the axioms and rules of classical propositional logic, let $\mathsf{S4}_{\sint}$ denote the $\mathsf{S4}$ axioms and rules for the $\sint$ modality, and let $\mathsf{S5}_{K}$ denote the $\mathsf{S5}$ axioms and rules for the $K$ modality \citep[see, e.g.,][]{FHMV}. Let $\mathbf{(KI)}$ denote the axiom scheme $K \phi \lthen \sint(\phi)$, and set
$$\mathsf{EL}_{\sint} \defeq \mathsf{CPL} + \mathsf{S4}_{\sint} + \mathsf{S5}_{K} + \mathbf{(KI)}.\footnote{This axiom system is closely related to an axiomatization presented by \citet{GP92} for a bimodal language containing both a ``local'' modality (quantifying over all accessible worlds in a relational structure) and a ``global'' modality (quantifying over \textit{all} worlds in the structure). In the subset space setting, the knowledge operator can be construed as a kind of global modality if one ignores the existence of states outside the current epistemic range. And indeed, the axiomatization of the universal modality given by \citeauthor{GP92} consists in the standard $\mathsf{S5}$ axioms together with what they call the ``inclusion'' axiom scheme, which corresponds exactly to our scheme $\mathbf{(KI)}$. Thanks to Ayb\"uke \"Ozg\"un for pointing out this connection.}$$

\begin{theorem} \label{thm:sac}
$\mathsf{EL}_{\sint}$ is a sound and complete axiomatization of $\mathcal{EL}_{\sint}$.
\end{theorem}

\begin{proof}
Soundness of $\mathsf{CPL} + \mathsf{S5}_{K}$ is easy to show in the usual way, while soundness of $\mathsf{S4}_{\sint} + \mathbf{(KI)}$ follows from Proposition \ref{pro:prp}.

Completeness can be proved by a relatively straightforward canonical model construction. Let $X$ denote the set of all maximal ($\mathsf{EL}_{\sint}$-)consistent subsets of $\mathcal{EL}_{\sint}$. Define a relation $\sim$ on $X$ by
$$x \sim y \; \dimp \; (\forall \phi \in \mathcal{EL}_{\sint})(K \phi \in x \dimp K \phi \in y).$$
Clearly $\sim$ is an equivalence relation; let $[x]$ denote the equivalence class of $x$ under $\sim$. These equivalence classes partition $X$ according to what is known, but we cannot simply take the set of epistemic ranges to be $\{[x] \: : \: x \in X\}$, since we require this set to be a topology on $X$ and to interact with the $\sint$ modality in the right way. So we need to do a bit more work to define $\T$.

For each $\phi \in \mathcal{EL}_{\sint}$, let $\widehat{\phi} \defeq \{x \in X \: : \: \phi \in x\}$. Roughly speaking, sets of the form $\widehat{\sint(\phi)}$ ought to be interiors in whatever topology we define; more precisely, if we have any hope of proving the Truth Lemma, below, then at a minimum we need to ensure that these sets are open. Thus, in order to respect both the $\sint$ and the $K$ modalities, we define
$$\B \defeq \big\{\widehat{\sint(\phi)} \cap [x] \: : \: \phi \in \mathcal{EL}_{\sint} \textrm{ and } x \in X\big\},$$
and let $\T$ be the topology generated by $\B$. In fact, it is not difficult to show (using $\mathsf{S4}_{\sint}$) that $\B$ is a basis for $\T$.\footnote{$\B$ is a basis for a topology $\T$ if every element of $\T$ is a union of elements of $\B$.}

For each $p \in \textsc{prop}$, set $v(p) \defeq \widehat{p}$. Let $\X = (X, \T, v)$. Clearly $\X$ is a topological subset model.

\begin{lemma}[Truth Lemma] \label{lem:tru}
For every $\phi \in \mathcal{EL}_{\sint}$, for all $x \in X$, $\phi \in x$ iff $(\X, x, [x]) \models \phi$.
\end{lemma}

\begin{proof}
First we note that $\widehat{\sint(\verum)} = X$, and thus for all $x \in X$ we have $[x] = \widehat{\sint(\verum)} \cap [x] \in \T$, so $(x, [x])$ is indeed an epistemic scenario of $\X$.

As usual, the proof proceeds by induction on the complexity of $\phi$. The base case holds by definition of $v$, and the inductive steps for the Boolean connectives are straightforward.

So suppose the result holds for $\phi$; let us show that it holds for $K \phi$. If $K \phi \in x$, then by definition of $\sim$ we know that $(\forall y \in [x])(K \phi \in y)$. But $K \phi \in y \rimp \phi \in y$, so $(\forall y \in [x])(\phi \in y)$, which by the inductive hypothesis implies that $(\forall y \in [x])((y,[y]) \models \phi)$. Since $[y] = [x]$, this is equivalent to $(\forall y \in [x])((y,[x]) \models \phi)$, which yields $(x,[x]) \models K \phi$.

For the converse, suppose that $K \phi \notin x$. Then $\{K \psi \: : \: K \psi \in x\} \cup \{\lnot \phi\}$ is consistent, for if not there is a finite subset $\Gamma \subseteq \{K \psi \: : \: K \psi \in x\}$ such that
$$\proves \bigwedge_{\chi \in \Gamma} \chi \lthen \phi$$
(where $\proves$ denotes provability in $\mathsf{EL}_{\sint}$), from which it follows (using $\textsf{S5}_{K}$) that
$$\proves \bigwedge_{\chi \in \Gamma} \chi \lthen K \phi,$$
which implies $K \phi \in x$, a contradiction. Therefore, we can extend $\{K \psi \: : \: K \psi \in x\} \cup \{\lnot \phi\}$ to some $y \in X$; by construction, we have $y \in [x]$ and $\phi \notin y$. This latter fact, by the inductive hypothesis, yields $(y,[y]) \not\models \phi$ and thus $(y,[x]) \not \models \phi$ (since $[x] = [y]$), whence $(x,[x]) \not\models K\phi$.

Now let us suppose that the result holds for $\phi$ and work to show that it also must hold for $\sint(\phi)$. If $\sint(\phi) \in x$, then observe that
$$x \in \widehat{\sint(\phi)} \cap [x] \subseteq \{y \in [x] \: : \: \phi \in y\};$$
this is an easy consequence of the fact that $\proves \sint(\phi) \lthen \phi$. Since $\widehat{\sint(\phi)} \cap [x]$ is open, it follows that
\begin{equation} \label{eqn:sin}
x \in \int(\{y \in [x] \: : \: \phi \in y\}).
\end{equation}
Now by the inductive hypothesis we have
\begin{eqnarray*}
\{y \in [x] \: : \: \phi \in y\} & = & \{y \in [x] \: : \: (y,[y]) \models \phi\}\\
& = & \{y \in [x] \: : \: (y,[x]) \models \phi\}\\
& = & \oval{\phi}^{[x]},
\end{eqnarray*}
which by (\ref{eqn:sin}) yields $x \in \int \oval{\phi}^{[x]}$, so $(x,[x]) \models \sint(\phi)$.

For the converse, suppose that $(x,[x]) \models \sint(\phi)$. Then $x \in \int \oval{\phi}^{[x]}$ which, as above, is equivalent to $x \in \int(\{y \in [x] \: : \: \phi \in y\})$. It follows that there is some basic open set $\widehat{\sint(\psi)} \cap [z]$ such that
$$x \in \widehat{\sint(\psi)} \cap [z] \subseteq \{y \in [x] \: : \: \phi \in y\};$$
of course, in this case it must be that $[z] = [x]$. This implies that for all $y \in [x]$, if $\sint(\psi) \in y$ then $\phi \in y$. From this we can deduce that
$$\{K \psi' \: : \: K \psi' \in x\} \cup \{\lnot(\sint(\psi) \lthen \phi)\}$$
is inconsistent, for if not it could be extended to a $y \in [x]$ with $\sint(\psi) \in y$ but $\phi \notin y$, a contradiction. Thus, we can find a finite subset $\Gamma \subseteq \{K \psi' \: : \: K \psi' \in x\}$ such that
$$\proves \bigwedge_{\chi \in \Gamma} \chi \lthen (\sint(\psi) \lthen \phi),$$
which implies (using $\mathsf{S5}_{K}$) that
$$\proves \bigwedge_{\chi \in \Gamma} \chi \lthen K(\sint(\psi) \lthen \phi).$$
This implies that $K(\sint(\psi) \lthen \phi) \in x$, so by $\mathbf{(KI)}$ we know also that $\sint(\sint(\psi) \lthen \phi) \in x$, from which it follows (using $\mathsf{S4}_{\sint}$) that $\sint(\psi) \lthen \sint(\phi) \in x$. Since $x \in \widehat{\sint(\psi)}$, we conclude that $\sint(\phi) \in x$, as desired.
\qed
\end{proof}

Completeness, of course, is an easy consequence: if $\phi$ is not a theorem of $\mathsf{EL}_{\sint}$, then $\{\lnot \phi\}$ is consistent and can be extended to some $x \in X$, in which case by Lemma \ref{lem:tru} we have $(\X,x,[x]) \not\models \phi$.
\qed
\end{proof}

Let $\mathsf{PAL}_{\sint}$ denote $\mathsf{EL}_{\sint}$ together with the six reduction schemes given in Proposition \ref{pro:red}.

\begin{corollary} \label{cor:sac}
$\mathsf{PAL}_{\sint}$ is a sound and complete axiomatization of $\mathcal{PAL}$ (with respect to $\sint$-semantics).
\end{corollary}

\begin{proof}
Soundness follows from soundness of $\mathsf{EL}_{\sint}$ together with Proposition \ref{pro:red}. For completeness, let $\phi$ be a valid $\mathcal{PAL}$ formula. Then we can find an $\mathcal{EL}_{\sint}$ formula $\tilde{\phi}$ such that $\proves_{\mathsf{PAL}_{\sint}} \phi \liff \tilde{\phi}$; this can be seen by running essentially the same argument presented in Proposition \ref{pro:eqv}, replacing $\models$ with $\proves_{\mathsf{PAL}_{\sint}}$. Now $\tilde{\phi}$ is valid because $\phi$ is, so by completeness of $\mathsf{EL}_{\sint}$ we can deduce that $\proves_{\mathsf{EL}_{\sint}} \tilde{\phi}$, and so $\proves_{\mathsf{PAL}_{\sint}} \tilde{\phi}$, hence $\proves_{\mathsf{PAL}_{\sint}} \phi$.
\qed
\end{proof}

\section{Discussion} \label{sec:dis}

Subset spaces are a natural setting in which to model the dynamics of knowledge. But the semantic tools they offer are not quite enough for a satisfying interpretation of public announcements. Intuitively, $\phi$ is \emph{announceable} exactly when \textit{some} state of knowledge entails $\phi$, but this notion of announceability need not itself be represented as a knowledge state, and so cannot in general serve as the foundation for an epistemic update.

Topological structure offers an elegant solution: the announceability of $\phi$ is realized as the \textit{topological interior} of (the extension of) $\phi$, which therefore becomes both the precondition for and the content of a successful announcement of $\phi$, as given by (\ref{eqn:ref}). Examples \ref{exa:tw3} and \ref{exa:jt2} show that this topological definition has significant advantages over the semantics proposed by \citet{WA13}: more formulas are announceable, and successful announcements have implications that go beyond the mere truth of the announced formula. Moreover, from a technical standpoint, a modalized interior operator is a familiar and well-studied object, so its central role in our semantics situates this work in the broad context of topological semantics for modal (and especially epistemic) logics.

The epistemic interpretation of the $\sint$ modality is of interest in this regard. In motivating our use of the interior operator, we touched on intuitions of ``knowability'', and indeed it is tempting to think of $\sint(\phi)$ as expressing that $\phi$ is knowable. However, this turns out to be problematic for essentially the same reasons that \emph{Moore formulas} are problematic \citep{Moore}. Recall the setting of \textit{The Target and the Wall} as depicted in Figure \ref{fgr:taw}, and consider the Moore formula $\mu = b \land \lnot K b$. It is easy to see that $(x,X) \models \sint(\mu)$ since, for instance, $x \in U \subseteq B = \oval{\mu}^{X}$. On the other hand, $(x,U) \not\models K \mu$; in fact, $K \mu$ entails both $Kb$ and $\lnot Kb$, a contradiction. So in this straightforward sense, $\mu$ is \textit{not} knowable.

Loosely speaking, this discrepancy stems from the kind of appeal being made to the state of knowledge $U$: although $U$ acts as a witness to $(x,X)$ satisfying $\sint(\mu)$, it is not, in this capacity, ever treated as the epistemic range with respect to which knowledge statements in the language are evaluated. As soon as it is, $\mu$ is falsified.

This distinction can be captured formally with an \emph{epistemic effort} modality as in the original development of subset space logic: \citet{DMP96} work with an enriched language including formulas of the form $\Diamond \phi$, interpreted by
\begin{equation*}
(\X,x,U) \models \Diamond \phi \quad \textrm{iff} \quad (\exists V \in \T)(x \in V \subseteq U \textrm{ and } (\X,x,V) \models \phi).
\end{equation*}
Such a formula might be read, ``after some (epistemic) effort, $\phi$ holds''. This makes the formula $\Diamond K \phi$ an intuitive candidate for expressing knowability, and the argument above demonstrates that $\sint(\mu)$ and $\Diamond K \mu$ are \textit{not} equivalent. Enriching our logical setting to include the effort modality would provide a formal framework in which to investigate the relationship between these two notions of knowability, and more generally between abstract epistemic effort and public announcements. This is the subject of ongoing research.

In a very similar vein, the link between knowability and announcements has been investigated by \citet{BBvDHHdL08}, who extend the syntax of the language of public announcements with an additional \emph{arbitrary announcement} modality we might denote by $[*]$; roughly speaking, $[*]\phi$ is true when all (suitably chosen) formulas $\psi$ are such that $[\psi]\phi$ holds. The dualized version $\langle * \rangle \phi$ is therefore naturally read as, ``there is an announcement after which $\phi$ is true''. This too yields a plausible candidate for knowability: $\langle * \rangle K \phi$, that which becomes known after some announcement \citep{vB04}. In recent work, \citet{vDKO14} extend the logical system we have developed here to include just such an arbitrary announcement modality (their work cites an earlier, unpublished draft of this paper \citep{Bjorndahl13}). 

Building on this work, \citet{vDKO15} extend the logic further to a multi-agent framework. Multi-agent extensions are valuable generalizations of any single-agent epistemic framework, but in this setting there may be a special significance for interpreting the $\sint$ modality. In our semantics, the epistemic range of a given epistemic scenario is keyed to the mental state of a particular agent---namely, the one who hears the announcement. But in the example above, the ``witness'' $U$ to $x$ being in $\int\oval{\mu}^{X}$ does not function as an epistemic range, but merely as an information set. This suggests that a more suitable reading for $\sint(\phi)$ might be, ``$\phi$ is knowable by some third party'', or perhaps even, ``$\phi$ is known to the one who made the announcement''. A multi-agent logic rich enough to represent public announcements \textit{along with their agential sources} (e.g., ``after an announcement of $\phi$ by agent $i$...'') might therefore be just the right setting in which to truly understand the epistemics of the $\sint$ modality.

\begin{acknowledgement}
The insight that subset spaces can be used to provide an elegant, model-internal mechanism for interpreting public announcements is due to \citet{WA13}. Their work inspired this paper. I am also indebted to Hans van Ditmarsch, Ayb\"uke \"Ozg\"un, and Kevin T.~Kelly for helpful discussion of this topic, and to Joseph Y.~Halpern for comments on an earlier draft.
\end{acknowledgement}

\bibliographystyle{plainnat}
\bibliography{../../Bibliography/abjorndahl}

\end{document}